\documentclass[a4paper,onecolumn,10pt]{amsart}

\usepackage{amsmath}
\usepackage{amssymb}  
\usepackage{amsthm}
\usepackage{amsfonts}
\usepackage{graphicx}
\usepackage{cancel}
\usepackage{bbm}
\usepackage{url}
\usepackage{enumerate} 
\usepackage{ upgreek }
\usepackage{dsfont}
\usepackage{color}
\usepackage{subcaption}

\usepackage{makecell}
\usepackage{diagbox}
\usepackage{wrapfig}
\usepackage{rotating}
\usepackage{tabularx}

\usepackage{hyperref}
\hypersetup{
    colorlinks=true,
    linkcolor=blue,
    citecolor=red,
    filecolor=magenta,      
    urlcolor=cyan,
}

\def\beq{\begin{equation}}
\def\eeq{\end{equation}}
\def\bq{\begin{quote}}
\def\eq{\end{quote}}
\def\ben{\begin{enumerate}}
\def\een{\end{enumerate}}
\def\bit{\begin{itemize}}
\def\eit{\end{itemize}}

\def\lset{\lbrace}
\def\rset{\rbrace}

\def\r|{\right|}

\newcommand\R{\mathbbm{R}}

\theoremstyle{plain}
\newtheorem{thm}{Theorem}

\newtheorem{defn}[thm]{Definition}

\theoremstyle{definition}

\begin{document}
\title{{Cutting cakes and kissing circles}}
\author{Alexander M\"uller-Hermes}
\address{\small{Institut Camille Jordan, Universit\'{e} Claude Bernard Lyon 1,\\ 43 boulevard du 11 novembre 1918, 69622 Villeurbanne cedex, France}}
\email{muellerh@posteo.net}

\date{\today}
\begin{abstract}
To divide a cake into equal sized pieces most people use a knife and a mixture of luck and dexterity. These attempts are often met with varying success. Through precise geometric constructions performed with the knife replacing Euclid's straightedge and without using a compass we find methods for solving certain cake-cutting problems exactly. Since it is impossible to exactly bisect a circular cake when its center is not known, our constructions need to use multiple cakes. Using three circular cakes we present a simple method for bisecting each of them or to find their centers. Moreover, given a cake with marked center we present methods to cut it into $n$ pieces of equal size for $n=3,4$ and $6$. Our methods are based upon constructions by Steiner and Cauer from the 19th and early 20th century.

\end{abstract}
\maketitle

\section{Introduction}

Dividing a cake among a given number of people is a relevant problem at many birthday parties. Mathematically, this problem has often been studied in the context of ``fair division'' going back to Hugo Steinhaus~\cite{steinhaus1949division}. Here, every cake-eating participant should receive a share they themselves consider to be fair. Many protocols for fair division have been studied (see for instance~\cite{dubins1961cut,stromquist1980cut,woodall1980dividing}), but they do not seem to be practical at birthday parties with more than two participants who might lack the enthusiasm to perform complicated protocols, or who are not able to accurately compare the sizes of different pieces of cake. For example, the author himself finds it even difficult to cut a cake into two pieces from which the larger piece cannot be selected immediately. In this article, we present precise geometric constructions for how to cut circular cakes into equal sized pieces. Our constructions are based on the following assumptions: 
\begin{enumerate}
\item Cakes are perfect circles.  
\item Straight lines can be carved into the cake (or the table) using a knife.
\end{enumerate}
It is our firm opinion that geometric considerations of cake-cutting should only use the knife, which is considered equivalent to a straightedge (i.e.~unmarked ruler) in Euclidean geometry. Non-circular cakes are beyond the scope of this article, but we challenge the reader to find techniques to cut whatever shape of bakery they might encounter into pieces of equal size.  

\section{History of cake cutting: The Poncelet-Steiner theorem}
\label{sec:PS}

With the assumptions stated above, we can cut cakes by geometric constructions following the Euclidean axioms, but without using the compass. To substitute for the compass the cakes themselves may be used as preexisting circles. The study of such constructions has a long history, and we will start with the following theorem proved by Jakob Steiner~\cite{steiner1833geometrischen} in 1833 after being conjectured by Victor Poncelet. We recommend~\cite{dorrie1965100great} for a well-written exposition of its proof.

\begin{thm}[Poncelet-Steiner]
Any construction in Euclidean geometry can be performed using the straightedge alone provided that a circle with its center is given.
\end{thm}

By the Poncelet-Steiner theorem it is possible to cut a cake with marked center into $n$ pieces of equal size whenever the regular $n$-gon is constructible in Euclidean geometry, i.e.~whenever $n$ is a product of distinct Fermat primes and a power of 2. 

Many cakes encountered in the wild have their center marked by some kind of decoration, but cakes with unmarked or inaccurately marked center are also quite common. Unfortunately, the conclusion of the Poncelet-Steiner theorem is no longer true, when the center of the given circle is not known. This observation is due to David Hilbert as noted in~\cite{cauer1912konstruktion}, and an exposition of his argument can be found in~\cite{rademacher1957enjoyment}. As a consequence, it is impossible to exactly cut a cake in half, when the center is not known and only a knife may be used. Luckily, having more cakes available saves the day as shown by Detlef Cauer~\cite{cauer1912konstruktion} in 1912. Focusing on situations most relevant for cutting cakes, it is possible to construct the center of a circle $c$ using the straightedge alone if another circle touching $c$ from the outside is given, or if two additional circles are given in any position such that neither of the three circles lies inside another. By the Poncelet-Steiner theorem this implies the following theorem.

\begin{thm}[Cauer]
Any construction in Euclidean geometry can be performed using the straightedge alone provided that either two circles are given touching from the outside, or three circles are given in any position such that neither of the three circles lies inside another. 
\end{thm}

\begin{figure}[hbt!]
        \center
        \includegraphics[scale=0.63]{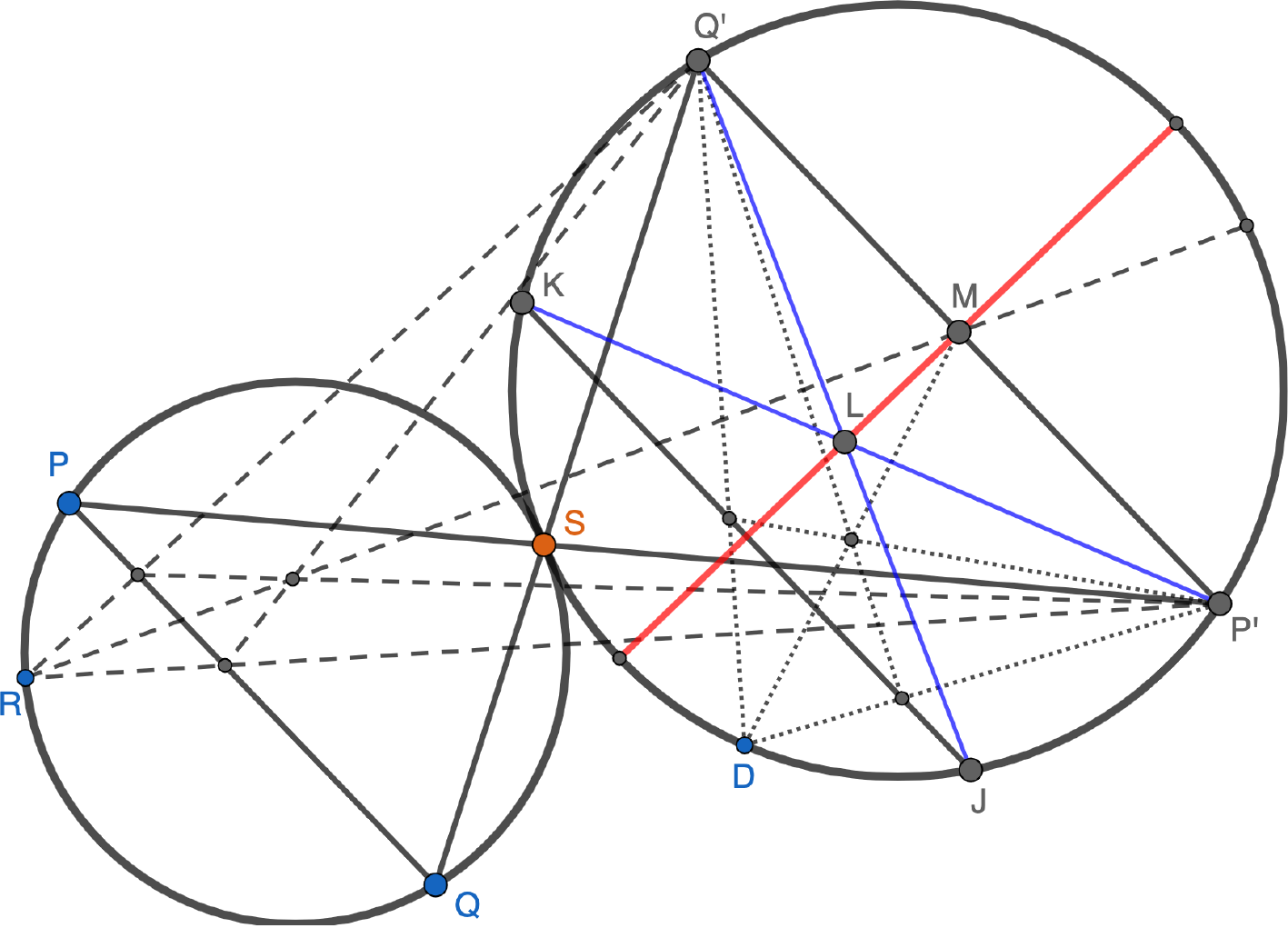}
        \caption{The red line indicates the final cut bisecting the larger cake.}
        \label{fig:CauerBisect}
\end{figure}

Since cakes, unlike circles drawn on a piece of paper, are movable objects, we can push two cakes together so that they touch in a point. By Cauer's theorem two cakes are then enough to perform general Euclidean constructions using a knife alone. Although we provide all ingredients for cutting pairs of touching cakes in this article, there is an important caveat: In the constructions known to the author (see for instance Figure \ref{fig:CauerBisect} for a bisection) auxiliary lines have to be carved into the table rather than into the surface of the cake. While we have not excluded this possibility in the assumptions stated above, it might create conflict in practice and the aspiring cake-cutter might not get invited anymore to parties with cakes to cut. In the following, we will present constructions that only use the surface of the cakes but require at least three cakes or a marked center. Specifically, we present such constructions for cutting a cake into $n$ pieces of equal size where $n=2,3,4$ and $6$. 

\section{Bisecting an odd number of cakes}
\label{sec:BisectOdd}

\begin{figure}[hbt!]
        \center
        \includegraphics[scale=0.6]{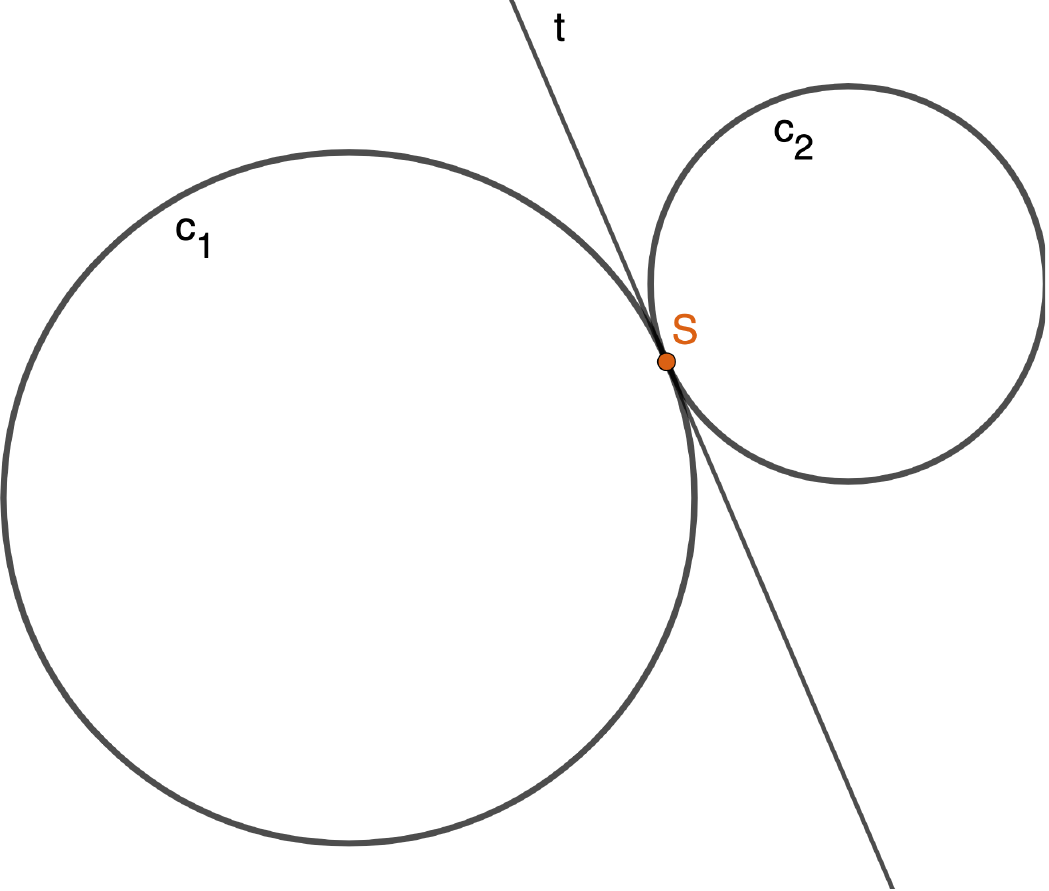}
        \caption{Kissing circles.}
        \label{fig:kissing}
\end{figure}

Our method for bisecting cakes is based on an elementary construction using kissing circles. We say that a pair of circles $c_1,c_2$ is \emph{kissing} in a point $S$ if the following conditions are satisfied:
\begin{enumerate}
\item The intersection satisfies $c_1\cap c_2=\lset S\rset$.
\item The circles $c_1$ and $c_2$ lie on opposite sides of their common tangent in $S$.
\end{enumerate}
See Figure \ref{fig:kissing} for an illustration.

Given a pair of kissing circles and a point $P$ on one of the circles we can construct a point $P'$ on the other circle by passing $P$ through the kissing point:

\begin{defn}[Passing through the kissing point]
Consider circles $c_1$ and $c_2$ kissing in a point $S$. For a point $P$ on $c_1$ we construct a point $P'$ on $c_2$ as follows:
\begin{enumerate}
\item If $P=S$, then we set $P'=S$.
\item If $P\neq S$, then we set $P'$ to be the intersection point different from $S$ of the line $PS$ with $c_{2}$. 
\end{enumerate}
We will say that $P'$ is obtained by \emph{passing $P$ through the kissing point} $S$.
\end{defn}

\begin{figure}[hbt!]
    \includegraphics[scale=0.06]{pic31.pdf}

    \includegraphics[scale=0.28]{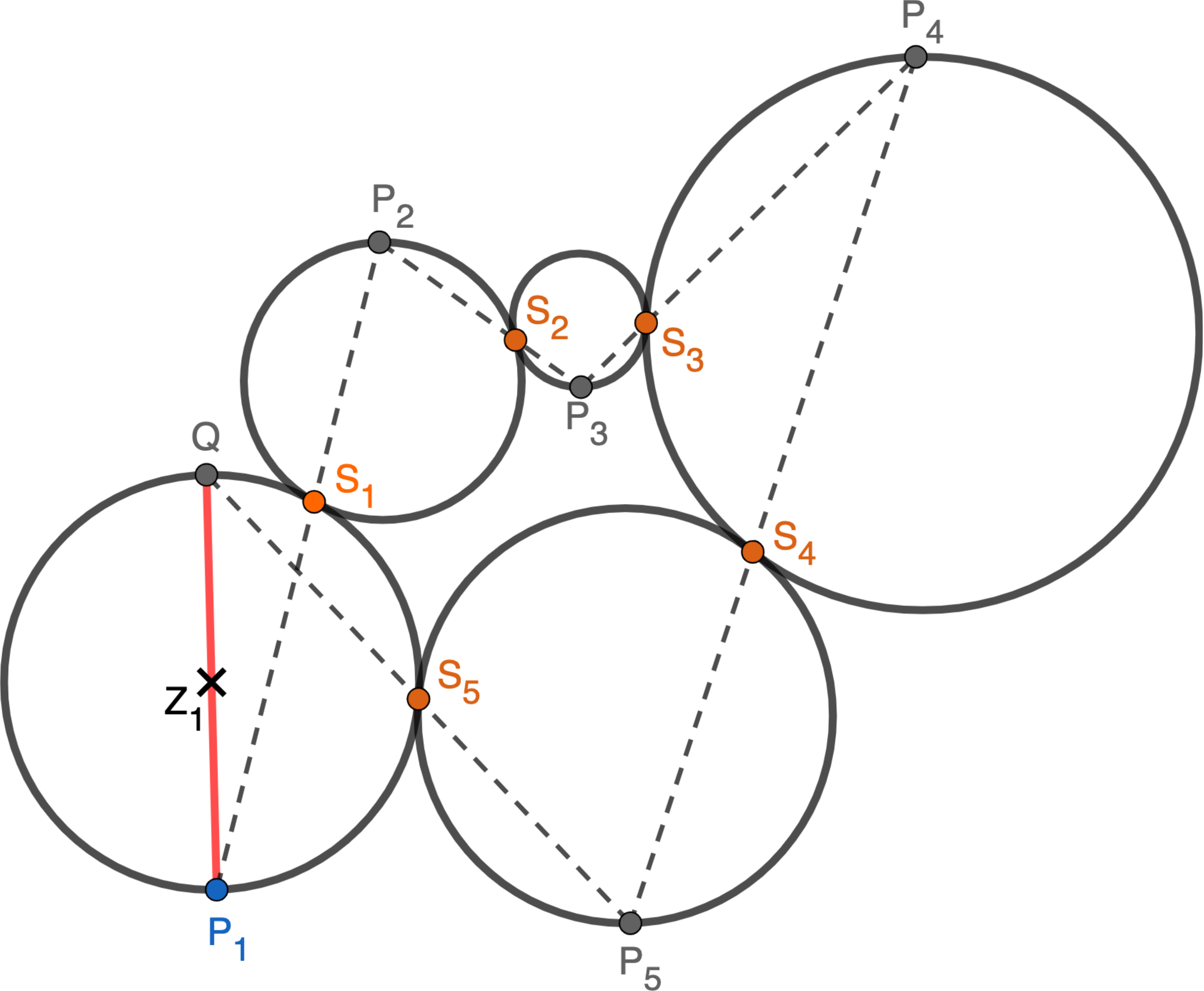}

   \caption{Construction of a diameter (red) using three and five circular cakes. Dotted lines should only be carved slightly into the surface of the cakes.}
   \label{fig:kissingCakesWithCenters}
\end{figure}

Suppose now that we are given an odd number of cakes. To bisect one of them we first push them together such that they form a closed chain of kissing circles as in Figure \ref{fig:kissingCakesWithCenters}, i.e.~circles $c_1,\ldots ,c_{2n+1}$ such that the pair $(c_{i}, c_{i+1})$ is kissing in the point $S_i$ and the pair $(c_{2n+1},c_1)$ is kissing in the point $S_{2n+1}$. Starting with a point $P_1$ on $c_1$, we construct a point $P_2$ on $c_2$ by passing $P_1$ through the kissing point $S_1$. Repeating this successively for the other circles leads to points $P_2,\ldots ,P_{2n+1}$ on $c_2,\ldots , c_{2n+1}$ as in Figure \ref{fig:kissingCakesWithCenters}. Finally, we construct a point $Q$ on $c_1$ by passing $P_{2n+1}$ through the kissing point $S_{2n+1}$. We claim that the chord $\overline{QP_{1}}$ is a diameter and hence bisects the first cake.

To make the previous statement precise we formulate the following theorem also covering the case of an even number of cakes. 

\begin{thm}
Consider circles $c_1,\ldots , c_n$ such that each pair $(c_i,c_{i+1})$ is kissing in a point $S_i$ and the pair $(c_n,c_{1})$ is kissing in a point $S_n$. Starting with a point $P_1$ on $c_1$ we construct points $P_2,\ldots , P_n$ such that each point $P_{i+1}$ is obtained by passing the point $P_i$ through the kissing point $S_i$. Let $Q$ be the point on $c_1$ obtained by passing $P_n$ through the kissing point $S_n$. Then, we have the following two statements:
\begin{enumerate}
\item If $n$ is even, then $Q=P_1$.
\item If $n$ is odd, then the chord $\overline{QP_1}$ is a diameter of $c_1$.
\end{enumerate} 
\label{thm:MainThm}
\end{thm}

While the construction outlined above works for every starting point $P_1$, we recommend to use $P_1=S_1$ for bisecting cakes in practice. By doing so, one needs to carve one auxiliary line less, and extending the final cut also bisects the second cake $c_2$. In the next section, we will present a proof of Theorem \ref{thm:MainThm}.

\section{Why does it work?}

To prove Theorem \ref{thm:MainThm}, we will apply the theory of dilations, i.e.~transformations of the plane transforming each line into a parallel line (where a line is considered parallel to itself). We will review all statements needed in the next paragraph, and we refer to~\cite[Section 13.2.]{coxeter1961introduction} and~\cite[Section 2.5]{aarts2009solid} for more details including proofs. 

There are two kinds of dilations: Translations of the plane, and central dilations (also known as homotheties) of the form $P\mapsto S + \kappa \overrightarrow{SP}$ with center $S$ and scale factor $\kappa\in \R$. The dilations form a group under composition. In particular, the composition of two central dilations with different centers $S_1$ and $S_2$ and scale factors $\kappa_1,\kappa_2$ is a translation when $\kappa_1\cdot\kappa_2=1$ or another central dilation with scale factor $\kappa_1\cdot\kappa_2\neq 1$. The composition of a translation and a central dilation with scale factor $\kappa$ is again a central dilation with the same scale factor $\kappa$ but different center, unless the translation is the identity. 

Given two circles $c_1$ and $c_2$ kissing in a point $S$, our operation ``passing through the kissing point'' introduced above and transforming $c_1$ into $c_2$ is a central dilation with center $S$ and a negative scale factor $\kappa=-r_2/r_1$, where $r_1$ and $r_2$ denote the radii of $c_1$ and $c_2$ respectively. Consider a closed chain of kissing circles $c_1,\ldots , c_n$ kissing in points $S_1,\ldots ,S_n$ as in Theorem \ref{thm:MainThm}, and let $f_1,\ldots ,f_n$ denote the central dilations such that $f_i$ passes points of $c_i$ through the kissing point $S_i$. Consider now the composition $g=f_n\circ \cdots \circ f_1$ of these central dilations. We know that $g$ is a dilation, and since $g(c_1)=c_1$ we conclude that $g$ is a central dilation with center $Z$ coinciding with the center of the circle $c_1$ and scale factor $\kappa$ either $+1$ or $-1$. Since each $f_i$ has a negative scaling factor, we see (by the composition rules outlined above) that $\kappa=-1$ when $n$ is odd. In this case the central dilation $g$ is the point reflection at the center $Z$ and $Q=g(P_1)$ is the point diametrically opposite of $P_1$ on $c_1$. On the other hand, if $n$ is even, then we see that $\kappa=+1$, and the central dilation $g$ is the identity and $Q=g(P_1)=P_1$. This proves Theorem \ref{thm:MainThm}. 

\section{Cutting cakes into more pieces}
\label{sec:CuttingMore}

So far, we have discussed how to find a diameter of a cake when there are at least three cakes available. Repeating the construction after rotating the cake, gives another diameter intersecting the first one in the center of the cake. When cutting three or more cakes we may therefore assume that the center is marked on each of them. By the Poncelet-Steiner theorem mentioned in Section \ref{sec:PS} it is then possible to cut each of the cakes into $3$, $4$ or $6$ pieces of equal size since the corresponding regular $n$-gons are constructible in Euclidean geometry. In the following, we will show how to actually do this in practice. 

Our constructions are based on two tricks developed by Jacob Steiner in~\cite{steiner1833geometrischen}, which can also be found in~\cite{dorrie1965100great}. The first trick constructs a line parallel to a given line $PQ$ through points $P$ and $Q$ when the midpoint $Z$ of the segment $\overline{PQ}$ is given. The second trick, closely related to the first, constructs the midpoint $Z$ of a segment $\overline{PQ}$ when a parallel line to the segment is given. Both statements follow directly from Ceva's theorem~\cite[p.220]{coxeter1961introduction} and we leave their proof to the reader.

\begin{figure}[hbt!]
		\center
        \includegraphics[scale=0.7]{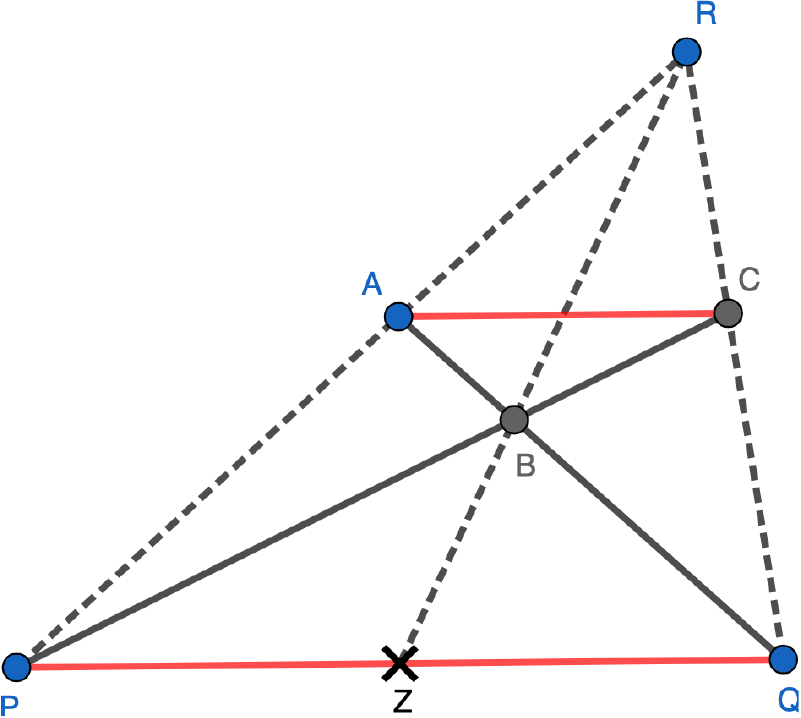}
        \caption{The red lines are parallel and $Z$ bisects the segment between $P$ and $Q$.}
        \label{fig:SteinerTrick}
\end{figure}

\begin{thm}[Two tricks by Steiner]
Consider points $P$ and $Q$. 
\begin{enumerate}
\item Let $Z$ be the point bisecting the segment $\overline{PQ}$ and consider a point $R$ not on the line $PQ$, and a point $A$ on the segment $\overline{PR}$. We construct a point $B$ as the intersection of $QA$ and $ZR$. Let $C$ denote the intersection of the line $PB$ with the line $QR$. Then, the line $AC$ is parallel to the line $PQ$.
\item Let points $A$ and $C$ be given such that the lines $AC$ and $PQ$ are parallel. We construct a point $R$ as the intersection of $PA$ and $QC$ and a point $B$ as the intersection of $PC$ and $QA$. Then, the line $RB$ intersects $PQ$ in the point $Z$ bisecting the segment $\overline{PQ}$
\end{enumerate} 
See Figure \ref{fig:SteinerTrick} for illustration.
\label{thm:SteinerTrick}
\end{thm}

\begin{figure}[hbt!]
		\center
        \includegraphics[scale=0.63]{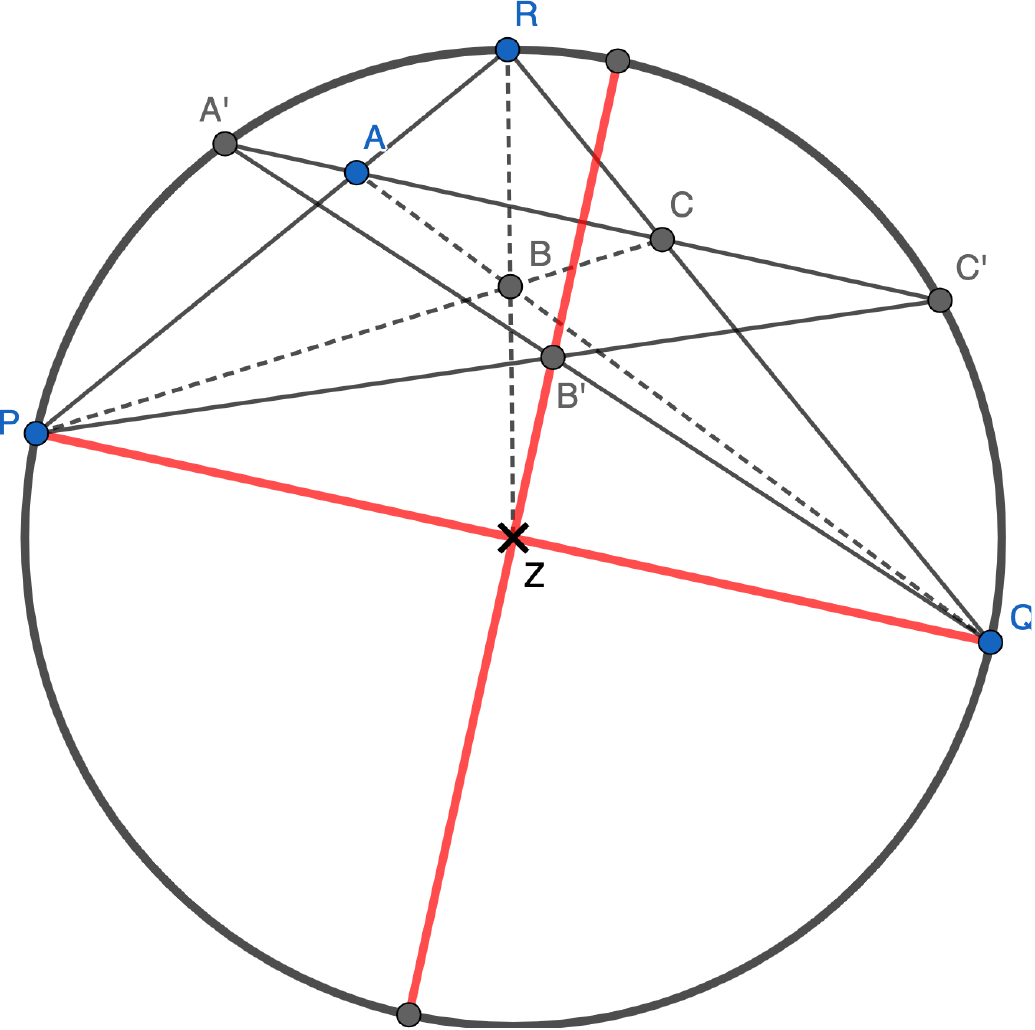}
        \caption{The red lines cut the cake into four pieces of equal size.}
        \label{fig:CuttinginFour}
\end{figure}

We will now show how to cut a cake with marked center $Z$ into $4$ pieces of equal size. Let $P$ and $Q$ denote a pair of antipodal points on the rim of the cake, i.e.~such that the chord $\overline{PQ}$ is a diameter. To cut the cake into four pieces of equal size we need to construct a diameter perpendicular to $\overline{PQ}$. Note that the center $Z$ bisects the diameter $\overline{PQ}$. By choosing a point $R$ on the rim of the cake and a point $A$ on the segment $\overline{PR}$, we can use the construction from the first case of Theorem \ref{thm:SteinerTrick} to find a point $C$ such that $AC$ is parallel to the $\overline{PQ}$ as in Figure \ref{fig:CuttinginFour}. 

Intersecting the line $AC$ with the circle yields points $A'$ and $C'$ such that the quadrilateral $QC'A'P$ is an isoceles trapezoid. Finally, the intersection point $B'$ of the diagonals $\overline{PC'}$ and $\overline{QA'}$ lies on a diameter perpendicular to $\overline{PQ}$ and by cutting along the lines $PQ$ and $ZB'$ we divide the cake into $4$ pieces of equal size.

\begin{figure}[hbt!]
		\center
        \includegraphics[scale=0.63]{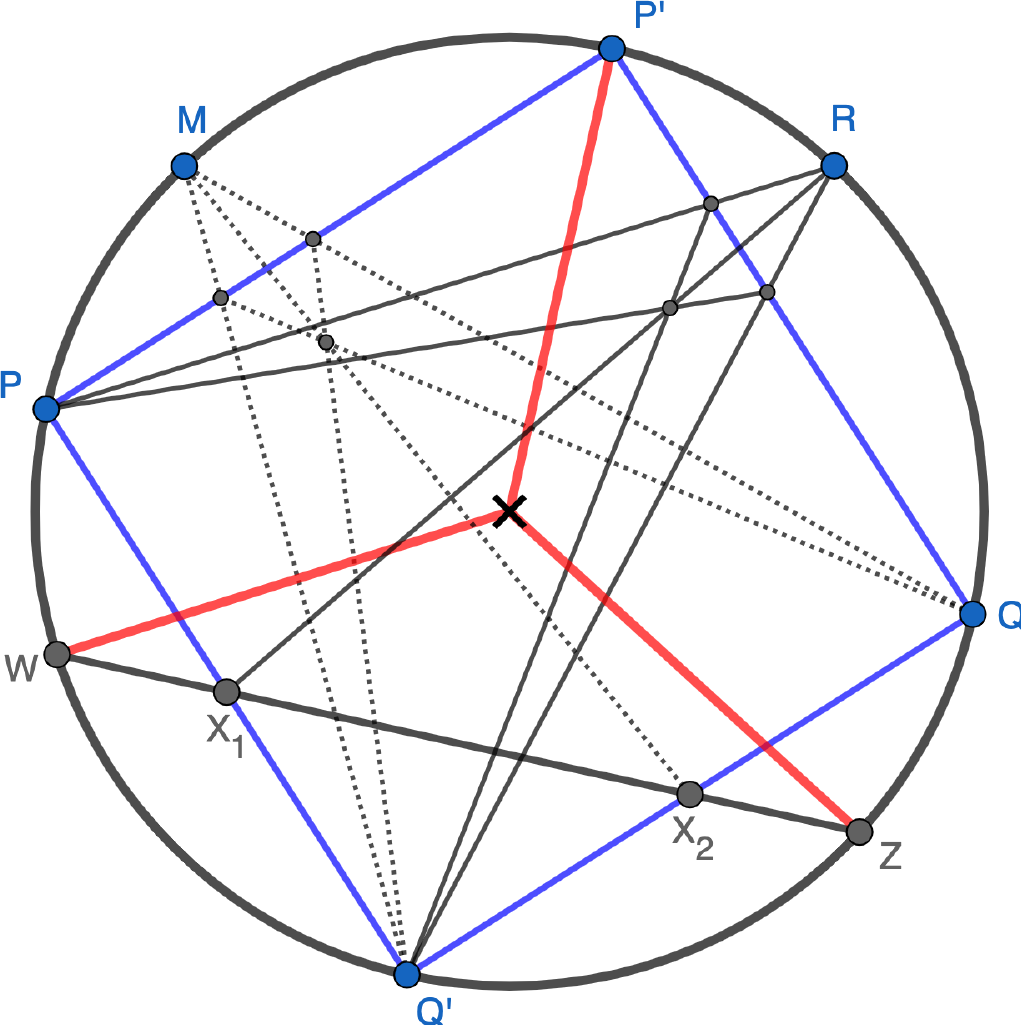}
        \caption{The red lines cut the cake into three pieces of equal size.}
        \label{fig:CuttinginThree}
\end{figure}

Next, we present a construction for cutting a cake into $3$ pieces of equal size. Again, we assume that a cake with marked center is given and using the previous construction we find four points $P,Q,P',Q'$ on its rim such that the quadrilateral $QP'PQ'$ is a square as in Figure \ref{fig:CuttinginThree}. We first note that the segments $\overline{PQ'}$ and $\overline{P'Q}$ are parallel. After choosing a point $R$ on the rim of the cake, we apply the construction from the second case of Theorem \ref{thm:SteinerTrick} to find a point $X_1$ bisecting the segment $\overline{PQ'}$. Repeating this construction for the parallel segments $\overline{QQ'}$ and $\overline{PP'}$ (and another point $M$ on the rim of the cake) we find the point $X_2$ bisecting the segment $\overline{QQ'}$. The line $X_1X_2$ intersects the rim in the points $W$ and $Z$ such that $P'WZ$ is an equilateral triangle. Finally, we cut from each point $P'$, $W$ and $Z$ to the center of the cake and obtain $3$ pieces of equal size. Note that the previous construction also allows to cut the cake into $6$ pieces of equal size by extending the red segments in Figure \ref{fig:CuttinginThree}.

\section*{Acknowledgements}

We thank Chris Perry and Emilie Elki\ae r for insightful comments and interesting discussions about cake-cutting that improved this article. Moreover, we thank Sergei Tabachnikov for pointing out a shorter proof of Theorem \ref{thm:MainThm} than the one stated in a previous version. We acknowledges financial support from the European Union's Horizon 2020 research and innovation programme under the Marie Sk\l odowska-Curie Action TIPTOP (grant no. 843414).

\bibliographystyle{alpha}
\bibliography{mybibliography.bib}

\begin{thebibliography}{Woo80}

\bibitem[Aar09]{aarts2009solid}
Jan~J. Aarts.
\newblock {\em Plane and Solid Geometry}.
\newblock Springer, 2009.

\bibitem[Cau12]{cauer1912konstruktion}
Detlef Cauer.
\newblock {{\"U}ber die Konstruktion des Mittelpunktes eines Kreises mit dem
  Lineal allein}.
\newblock {\em Mathematische Annalen}, 73(1):90--94, 1912.

\bibitem[Cox69]{coxeter1961introduction}
Harold S.~M. Coxeter.
\newblock {\em Introduction to geometry}.
\newblock {John Wiley \& Sons, Inc.}, 1969.

\bibitem[D{\"o}r65]{dorrie1965100great}
Heinrich D{\"o}rrie.
\newblock {\em {100 Great Problems of Elementary Mathematics: Their History and
  Solution}}.
\newblock Dover, 1965.

\bibitem[DS61]{dubins1961cut}
Lester~E. Dubins and Edwin~H. Spanier.
\newblock How to cut a cake fairly.
\newblock {\em The American Mathematical Monthly}, 68(1P1):1--17, 1961.

\bibitem[RT57]{rademacher1957enjoyment}
Hans Rademacher and Otto Toeplitz.
\newblock {\em {Enjoyment of Mathematics: Selections from Mathematics for the
  Amateur}}.
\newblock Princeton University Press, 1957.

\bibitem[Ste33]{steiner1833geometrischen}
Jacob Steiner.
\newblock {\em {Die geometrischen Konstructionen ausgef\"uhrt mittelst der
  geraden Linie und eines festen Kreises}}, volume~60.
\newblock D{\"u}mmler, 1833.

\bibitem[Ste49]{steinhaus1949division}
Hugo Steinhaus.
\newblock Sur la division pragmatique.
\newblock {\em Econometrica: Journal of the Econometric Society}, pages
  315--319, 1949.

\bibitem[Str80]{stromquist1980cut}
Walter Stromquist.
\newblock How to cut a cake fairly.
\newblock {\em The American Mathematical Monthly}, 87(8):640--644, 1980.

\bibitem[Woo80]{woodall1980dividing}
Douglas~R. Woodall.
\newblock Dividing a cake fairly.
\newblock {\em Journal of Mathematical Analysis and Applications},
  78(1):233--247, 1980.

\end{thebibliography}

\end{document}